# ESTIMATION OF NONLINEAR MODELS WITH BERKSON MEASUREMENT ERRORS[1]

By Liqun Wang

*University of Manitoba*

This paper is concerned with general nonlinear regression models where the predictor variables are subject to Berkson-type measurement errors. The measurement errors are assumed to have a general parametric distribution, which is not necessarily normal. In addition, the distribution of the random error in the regression equation is nonparametric. A minimum distance estimator is proposed, which is based on the first two conditional moments of the response variable given the observed predictor variables. To overcome the possible computational difficulty of minimizing an objective function which involves multiple integrals, a simulation-based estimator is constructed. Consistency and asymptotic normality for both estimators are derived under fairly general regularity conditions.

**1. Introduction.** In many scientific studies researchers are interested in the nonlinear relationship

$$Y = g(X;\theta) + \varepsilon, \tag{1}$$

where $Y \in \mathbb{R}$ is the response variable, $X \in \mathbb{R}^k$ is the predictor variable, $\theta \in \mathbb{R}^p$ is the unknown regression parameter and $\varepsilon$ is the random error. In many experiments, it is too costly or impossible to measure the predictor $X$ exactly. Instead, a proxy $Z$ of $X$ is measured.

For example, an epidemiologist studies the severity of a lung disease, $Y$, among the residents in a city in relation to the amount of certain air pollutants, $X$. Assume the air pollutants are measured at certain observation stations in the city. The actual exposure of the residents to the pollutants $X$, however, may vary randomly from the values $Z$ measured at these stations.

Received March 2003; revised August 2003.
[1]Supported by the Natural Sciences and Engineering Research Council of Canada.

In this case, $X$ can be expressed as $Z$ plus a random error, which represents the individual variation in the exposure from the measured exposure.

Other examples include agricultural or medical studies, where the relations between the yield of a crop or the efficacy of a drug, $Y$, and the amount of a fertilizer or drug used, $X$, are studied. Suppose the fertilizer or the drug is applied at predetermined doses $Z$. The actual absorption of the fertilizer in the crop or the drug in the patients' blood, however, may vary randomly around the set doses, because of the local earth conditions or the individual biological conditions. In these cases, if the amount $Z$ is properly calibrated, then the actual absorption $X$ will vary around $Z$ randomly, so that in average the random variation $X - Z$ will be zero.

In all situations mentioned above, a reasonable model for the measurement errors is the so-called Berkson model

$$X = Z + \delta, \tag{2}$$

where $\delta$ is the unobserved random measurement error which is assumed to be independent of the observed predictor variable $Z$. More explanations and motivations of the Berkson-error model can be found in Fuller [(1987), pages 79 and 80].

The stochastic structure of the Berkson measurement error model (2) is fundamentally different from the classical errors-in-variables model, where the measurement error is independent of $X$, but dependent on $Z$. This distinctive feature leads to completely different procedures in parameter estimation and inference for the models.

Estimation of the linear Berkson measurement error models is discussed in Fuller [(1987), pages 81–83] and Cheng and Van Ness [(1995), pages 35–38]. For nonlinear models, an approximative method called regression calibration is presented by Carroll, Ruppert and Stefanski [(1995), Chapter 3]. Recently, Huwang and Huang (2000) studied a univariate polynomial model where $g(x;\theta)$ is a polynomial in $x$ of a known order and showed that the least squares estimators based on the first two conditional moments of $Y$ given $Z$ are consistent. Wang (2003) considered general univariate nonlinear models where all random errors are normally distributed and showed that the minimum distance estimator based on the first two conditional moments of $Y$ given $Z$ is consistent and asymptotically normally distributed.

In many practical applications, however, there is often more than one predictor variable which is subject to measurement errors. Moreover, the random errors $\varepsilon$ and $\delta$ may have distributions other than the normal distribution. The goal of this paper is to generalize the results of Wang (2003) to the nonlinear models with multivariate predictor variables, where the measurement error $\delta$ has a general parametric distribution $f_\delta(t;\psi)$, $\psi \in \Psi \subset \mathbb{R}^q$, and the random error $\varepsilon$ has a nonparametric distribution with mean zero and variance $\sigma_\varepsilon^2$. Thus, (1) and (2) represent a semiparametric model. Our main interest is to estimate parameters $\gamma = (\theta', \psi', \sigma_\varepsilon^2)'$. We show that the minimum distance estimator of Wang (2003) is still consistent and asymptotically normally distributed. For the general model in this paper, however, a computational issue arises, because the objective function to be minimized involves multiple integrals for which explicit forms may not always be obtained. To overcome this difficulty, we propose a simulation-based estimator which is shown to be consistent and asymptotically normally distributed under regularity conditions similar to those for the minimum distance estimator.

Throughout the paper we assume that $Z, \delta$ and $\varepsilon$ are independent and $Y$ has finite second moment. In addition, we adopt the common assumption in the literature that the measurement error is "nondifferential" in the sense that the conditional expectation of $Y$ given $X$ and $Z$ is the same as the conditional expectation of $Y$ given $X$. Although in this paper $Z$ is assumed to be a random variable, it is easy to see that all results continue to hold if the observations of $Z$, $Z_1, Z_2, \dots, Z_n$, are treated as fixed constants such that the limits $\lim_{n\to\infty} \sum_{i=1}^n Z_i/n$ and $\lim_{n\to\infty} \sum_{i=1}^n Z_i Z_i'/n$ are finite.

The paper is organized as follows. In Section 2 we give three examples to motivate our estimation method. In Section 3 we formally define the minimum distance estimator and derive its consistency and asymptotic normality under some regularity conditions. In Section 4 we propose a simulation-based estimator and derive its consistency and asymptotic normality. Finally, conclusion and discussion are given in Section 5, whereas proofs of the theorems are given in Section 6.

## 2. Examples.

To motivate our estimation method, let us consider some examples. To simplify notation, let us consider the case where the measurement error $\delta = (\delta_1, \delta_2, \dots, \delta_k)'$ has the normal distribution $N(0, \sigma_\delta^2 I_k)$, where $0 < \sigma_\delta^2 < \infty$ and $I_k$ is the $k$-dimensional identity matrix.

EXAMPLE 1. First consider the model $g(x;\theta) = \theta_1 x_1 + \theta_3 e^{\theta_2 x_2}$, where $\theta_2\theta_3 \neq 0$. For this model the conditional moment of $Y$ given $Z$ can be written as

$$
\begin{aligned}
E(Y|Z) &= \theta_1 Z_1 + \theta_1 E(\delta_1) + \theta_3 e^{\theta_2 Z_2} E(e^{\theta_2\delta_2}) \\
&= \varphi_1 Z_1 + \varphi_3 e^{\varphi_2 Z_2},
\end{aligned}
\tag{3}
$$

where $\varphi_1 = \theta_1$, $\varphi_2 = \theta_2$ and $\varphi_3 = \theta_3 e^{\theta_2^2\sigma_\delta^2/2}$. Similarly, the second conditional moment of $Y$ given $Z$ can be written as

$$
\begin{aligned}
E(Y^2|Z) &= \theta_1^2 E[(Z_1+\delta_1)^2|Z] + \theta_3^2 E[e^{2\theta_2(Z_2+\delta_2)}|Z] \\
&\quad + 2\theta_1\theta_3 E[(Z_1+\delta_1)e^{\theta_2(Z_2+\delta_2)}|Z] + E(\varepsilon^2)
\end{aligned}
$$

$$
\begin{aligned}
&= \theta_1^2(Z_1^2+\sigma_\delta^2)+\theta_3^2 e^{2\theta_2 Z_2}E(e^{2\theta_2\delta_2}) \\
&\quad + 2\theta_1\theta_3 Z_1 e^{\theta_2 Z_2}E(e^{\theta_2\delta_2})+\sigma_\varepsilon^2 \\
&= \varphi_4+\varphi_1^2 Z_1^2+\varphi_5 e^{2\varphi_2 Z_2}+2\varphi_1\varphi_3 Z_1 e^{\varphi_2 Z_2},
\end{aligned} \tag{4}
$$

where $\varphi_4 = \theta_1^2\sigma_\delta^2 + \sigma_\varepsilon^2$ and $\varphi_5 = \theta_3^2 e^{2\theta_2^2\sigma_\delta^2}$. Since (3) and (4) are the usual nonlinear regression equations and both $Y$ and $Z$ are observable, $(\varphi_i)$ are identified by these equations and, therefore, can be consistently estimated using the nonlinear least squares method. Furthermore, the original parameters $(\theta_i, \sigma_\delta^2, \sigma_\varepsilon^2)$ are identified because the mapping $(\theta_i, \sigma_\delta^2, \sigma_\varepsilon^2) \mapsto (\varphi_i)$ is bijective. Indeed, it is straightforward to calculate that $\theta_1 = \varphi_1$, $\theta_2 = \varphi_2$, $\theta_3 = \varphi_3^2/\sqrt{\varphi_5}$, $\sigma_\delta^2 = \log(\varphi_5/\varphi_3^2)/\varphi_2^2$ and $\sigma_\varepsilon^2 = \varphi_4 - \varphi_1^2\log(\varphi_5/\varphi_3^2)/\varphi_2^2$.

EXAMPLE 2. Now consider another model $g(x;\theta) = \theta_1\exp(x'\theta_2)$, where $\theta_1 \neq 0$, $0 \neq \theta_2 \in \mathbb{R}^{p-1}$ and $p > 1$. For this model the first conditional moment of $Y$ given $Z$ can be written as

$$
\begin{aligned}
E(Y|Z) &= \theta_1 e^{Z'\theta_2}E(e^{\delta'\theta_2}) \\
&= \varphi_1 e^{Z'\varphi_2},
\end{aligned} \tag{5}
$$

where $\varphi_1 = \theta_1\exp(\theta_2'\theta_2\sigma_\delta^2/2)$ and $\varphi_2 = \theta_2$. The second conditional moment is given by

$$
\begin{aligned}
E(Y^2|Z) &= \theta_1^2 e^{2Z'\theta_2}E(e^{2\delta'\theta_2}) + E(\varepsilon^2) \\
&= \varphi_3 e^{2Z'\varphi_2} + \varphi_4,
\end{aligned} \tag{6}
$$

where $\varphi_3 = \theta_1^2 e^{2\theta_2'\theta_2\sigma_\delta^2}$ and $\varphi_4 = \sigma_\varepsilon^2$. Again, $(\varphi_i)$ are identified by (5) and (6) and the nonlinear least squares method. Furthermore, the original parameters $(\theta_i, \sigma_\delta^2, \sigma_\varepsilon^2)$ are identified because the mapping $(\theta_i, \sigma_\delta^2, \sigma_\varepsilon^2) \mapsto (\varphi_i)$ is bijective. Indeed, straightforward calculation shows that $\theta_1 = \varphi_1^2/\sqrt{\varphi_3}$, $\theta_2 = \varphi_2$, $\sigma_\delta^2 = \log(\varphi_3/\varphi_1^2)/\varphi_2'\varphi_2$ and $\sigma_\varepsilon^2 = \varphi_4$.

EXAMPLE 3. Further, let us consider the polynomial model $g(x;\theta) = \theta_1 x_1 + \theta_2 x_2 + \theta_3 x_1^2 + \theta_4 x_2^2 + \theta_5 x_1 x_2$. For this model the first two conditional moments are, respectively,

$$
E(Y|Z) = (\theta_3+\theta_4)\sigma_\delta^2 + \theta_1 Z_1 + \theta_2 Z_2 + \theta_3 Z_1^2 + \theta_4 Z_2^2 + \theta_5 Z_1 Z_2 \tag{7}
$$

and

$$
\begin{aligned}
E(Y^2|Z) &= E[g^2(Z+\delta;\theta)|Z] + E(\varepsilon^2|Z) \\
&= E[g^2(Z+\delta;\theta)|Z] + \sigma_\varepsilon^2.
\end{aligned} \tag{8}
$$

Again, it is easy to see that parameters $(\theta_i)$ and $\sigma_\delta^2$ are identified by the nonlinear regression equation (7), whereas $\sigma_\varepsilon^2$ is identified by (8). Thus, all parameters in this model can be consistently estimated using the first two conditional moment equations.

The above examples suggest that in many situations, parameters in nonlinear models can be identified and, therefore, consistently estimated using the first two conditional moments of $Y$ given $Z$. The fact that parameters of Berkson measurement error models may be identified in nonlinear regression was first noticed by Rudemo, Ruppert and Streibig (1989). The identifiability of the univariate polynomial model was shown by Huwang and Huang (2000). For general nonlinear models, it is worthwhile noting that even in the case where the mapping $(\theta_i, \sigma_\delta^2, \sigma_\varepsilon^2) \mapsto (\varphi_i)$ is not bijective, the original parameters $(\theta_i, \sigma_\delta^2, \sigma_\varepsilon^2)$ can still be identified, if appropriate restrictions on them are imposed. In the next section, we develop a minimum distance estimator for the general nonlinear model (1) and (2) based on the first two conditional moments and derive its asymptotic properties.

**3. Minimum distance estimator.** Under the assumptions for model (1) and (2), the first two conditional moments of $Y$ given $Z$ are respectively given by

$$\begin{aligned} E(Y|Z) &= E[g(Z+\delta;\theta)|Z] + E(\varepsilon|Z) \\ &= \int g(Z+t;\theta) f_\delta(t;\psi)\,dt \end{aligned} \tag{9}$$

and

$$\begin{aligned} E(Y^2|Z) &= E[g^2(Z+\delta;\theta)|Z] + E(\varepsilon^2|Z) \\ &= \int g^2(Z+t;\theta) f_\delta(t;\psi)\,dt + \sigma_\varepsilon^2. \end{aligned} \tag{10}$$

Throughout this paper, unless otherwise stated explicitly, all integrals are taken to be over the space $\mathbb{R}^k$. Further, let $\gamma = (\theta', \psi', \sigma_\varepsilon^2)'$ denote the vector of model parameters and let $\Gamma = \Theta \times \Psi \times \Sigma \subset \mathbb{R}^{p+q+1}$ denote the parameter space. The true parameter value of model (1) and (2) is denoted by $\gamma_0 \in \Gamma$. For every $z \in \mathbb{R}^k$ and $\gamma \in \Gamma$, define

$$m_1(z;\gamma) = \int g(z+t;\theta) f_\delta(t;\psi)\,dt, \tag{11}$$

$$m_2(z;\gamma) = \int g^2(z+t;\theta) f_\delta(t;\psi)\,dt + \sigma_\varepsilon^2. \tag{12}$$

Then $m_1(Z;\gamma_0) = E(Y|Z)$ and $m_2(Z;\gamma_0) = E(Y^2|Z)$.

Now suppose $(Y_i, Z_i')', i = 1, 2, \ldots, n$, is an i.i.d. random sample and let

$$\rho(Y_i, Z_i;\gamma) = (Y_i - m_1(Z_i;\gamma), Y_i^2 - m_2(Z_i;\gamma))'.$$

Then the minimum distance estimator (MDE) $\hat{\gamma}_n$ for $\gamma$ based on moment equations (9) and (10) is defined as

$$\hat{\gamma}_n = \underset{\gamma\in\Gamma}{\arg\min}\, Q_n(\gamma),$$

where the objective function

$$Q_n(\gamma) = \sum_{i=1}^{n} \rho'(Y_i, Z_i; \gamma) W(Z_i) \rho(Y_i, Z_i; \gamma) \tag{13}$$

and $W(Z_i)$ is a $2 \times 2$ weighting matrix which may depend on $Z_i$.

Regularity conditions under which $\hat{\gamma}_n$ is identified, consistent and asymptotically normally distributed are well known in the nonlinear regression literature; see, for example, Amemiya [(1985), Chapter 5], Gallant [(1987), Chapter 5] and Seber and Wild [(1989), Chapter 12]. Usually these conditions are expressed in a variety of forms.

In the following, we adopt the setup of Amemiya (1985) and express these regularity conditions in terms of the regression function $g(x;\theta)$ and measurement error distribution $f_\delta(t;\psi)$. Let $\mu$ denote the Lebesgue measure and let $\|\cdot\|$ denote the Euclidean norm in $\mathbb{R}^d$. Then we assume the following conditions for the consistency of the MDE $\hat{\gamma}_n$.

ASSUMPTION A1. $g(x;\theta)$ is a measurable function of $x$ for every $\theta \in \Theta$, and is continuous in $\theta \in \Theta$, for $\mu$-almost all $x$.

ASSUMPTION A2. $f_\delta(t;\psi)$ is continuous in $\psi \in \Psi$ for $\mu$-almost all $t$.

ASSUMPTION A3. The parameter space $\Gamma \subset \mathbb{R}^{p+q+1}$ is compact.

ASSUMPTION A4. The weight $W(Z)$ is nonnegative definite with probability 1 and satisfies $E\|W(Z)\| < \infty$.

ASSUMPTION A5. $\int \sup_\Psi f_\delta(t;\psi)\,dt < \infty$ and $E(\|W(Z)\|+1) \times \int \sup_{\Theta\times\Psi} g^4(Z+t;\theta) f_\delta(t;\psi)\,dt < \infty$.

ASSUMPTION A6. $E[\rho(Y,Z;\gamma) - \rho(Y,Z;\gamma_0)]' W(Z)[\rho(Y,Z;\gamma) - \rho(Y,Z;\gamma_0)] = 0$ if and only if $\gamma = \gamma_0$.

The above regularity conditions are common in the literature of nonlinear regression. In particular, Assumptions A1 and A2 are usually used to ensure that the objective function $Q_n(\gamma)$ is continuous in $\gamma$. Similarly, the compactness of the parameter space $\Gamma$ is often assumed. From a practical point of view, Assumption A3 is not as restrictive as it seems to be, because for any given problem one usually has some information about the possible range of the parameters. Assumption A5 contains moment conditions which imply the uniform convergence of $Q_n(\gamma)$. In view of (9) and (10), this assumption means that $Y$ and $\varepsilon$ have finite fourth moments. It is easy to see that Assumptions A1, A2 and A5 are satisfied, if $g(x;\theta)$ is a polynomial in $x$ and the measurement error $\delta$ has a normal distribution. Finally, Assumption A6 is the usual condition for identifiability of parameters, which means that the true parameter value $\gamma_0$ is the unique minimizer of the objective function $Q_n(\gamma)$ for large $n$.

THEOREM 1. *Under Assumptions A1–A6, the MDE* $\hat{\gamma}_n \xrightarrow{P} \gamma_0$, *as* $n \to \infty$.

To derive the asymptotic distribution for the MDE $\hat{\gamma}_n$, we assume further regularity conditions as follows.

ASSUMPTION A7. There exist open subsets $\theta_0 \in \Theta_0 \subset \Theta$ and $\psi_0 \in \Psi_0 \subset \Psi$, in which $g(x;\theta)$ is twice continuously differentiable with respect to $\theta$ and $f_\delta(t;\psi)$ is twice continuously differentiable with respect to $\psi$, for $\mu$-almost all $x$ and $t$, respectively. Furthermore, their first two derivatives satisfy

$$E\|W(Z)\| \int \sup_{\Theta_0\times\Psi_0} \left\| \frac{\partial g(Z+t;\theta)}{\partial\theta} \right\|^2 f_\delta(t;\psi)\,dt < \infty,$$

$$E\|W(Z)\| \int \sup_{\Theta_0\times\Psi_0} \left\| \frac{\partial^2 g(Z+t;\theta)}{\partial\theta\,\partial\theta'} \right\|^2 f_\delta(t;\psi)\,dt < \infty,$$

$$E\|W(Z)\| \int \sup_{\Psi_0} \left\| \frac{\partial f_\delta(t;\psi)}{\partial\psi} \right\| dt < \infty,$$

$$E\|W(Z)\| \int \sup_{\Theta_0\times\Psi_0} g^2(Z+t;\theta) \left\| \frac{\partial f_\delta(t;\psi)}{\partial\psi} \right\| dt < \infty,$$

$$E\|W(Z)\| \int \sup_{\Theta_0\times\Psi_0} \left\| \frac{\partial g(Z+t;\theta)}{\partial\theta} \right\|^2 \left\| \frac{\partial f_\delta(t;\psi)}{\partial\psi} \right\| dt < \infty,$$

$$E\|W(Z)\| \int \sup_{\Theta_0\times\Psi_0} |g(Z+t;\theta)| \left\| \frac{\partial^2 f_\delta(t;\psi)}{\partial\psi\,\partial\psi'} \right\| dt < \infty,$$

$$E\|W(Z)\| \int \sup_{\Theta_0\times\Psi_0} g^2(Z+t;\theta) \left\| \frac{\partial^2 f_\delta(t;\psi)}{\partial\psi\,\partial\psi'} \right\| dt < \infty.$$

ASSUMPTION A8. The matrix

$$B = E\left[ \frac{\partial\rho'(Y,Z;\gamma_0)}{\partial\gamma} W(Z) \frac{\partial\rho(Y,Z;\gamma_0)}{\partial\gamma'} \right]$$

is nonsingular, where

$$\frac{\partial\rho'(Y,Z;\gamma_0)}{\partial\gamma} = -\left( \frac{\partial m_1(Z;\gamma_0)}{\partial\gamma}, \frac{\partial m_2(Z;\gamma_0)}{\partial\gamma} \right).$$

Again, Assumptions A7 and A8 are commonly seen regularity conditions which are sufficient for the asymptotic normality of the minimum distance estimators. Assumption A7 ensures that the first derivative of $Q_n(\gamma)$ admits a first-order Taylor expansion and the second derivative of $Q_n(\gamma)$ converges uniformly. This assumption and the dominated convergence theorem together imply that the first derivatives $\partial m_1(z;\gamma)/\partial\gamma$ and $\partial m_2(z;\gamma)/\partial\gamma$ exist and their elements are respectively given by

$$\frac{\partial m_1(z;\gamma)}{\partial\theta} = \int \frac{\partial g(z+t;\theta)}{\partial\theta} f_\delta(t;\psi)\,dt,$$
$$\frac{\partial m_1(z;\gamma)}{\partial\psi} = \int g(z+t;\theta)\frac{\partial f_\delta(t;\psi)}{\partial\psi}\,dt,$$
$$\frac{\partial m_1(z;\gamma)}{\partial\sigma_\varepsilon^2} = 0$$

and

$$\frac{\partial m_2(z;\gamma)}{\partial\theta} = 2\int \frac{\partial g(z+t;\theta)}{\partial\theta} g(z+t;\theta) f_\delta(t;\psi)\,dt,$$
$$\frac{\partial m_2(z;\gamma)}{\partial\psi} = \int g^2(z+t;\theta)\frac{\partial f_\delta(t;\psi)}{\partial\psi}\,dt,$$
$$\frac{\partial m_2(z;\gamma)}{\partial\sigma_\varepsilon^2} = 1.$$

Finally, Assumption A8 implies that the second derivative of $Q_n(\gamma)$ has a nonsingular limiting matrix. Again, Assumptions A7 and A8 are satisfied for the polynomial model $g(x;\theta)$ and the normal measurement error $\delta$.

THEOREM 2. *Under Assumptions A1–A8, as* $n\to\infty$, $\sqrt{n}(\hat\gamma_n-\gamma_0)\xrightarrow{L} N(0,B^{-1}CB^{-1})$, *where*

$$C = E\left[\frac{\partial\rho'(Y,Z;\gamma_0)}{\partial\gamma}W(Z)\rho(Y,Z;\gamma_0)\rho'(Y,Z;\gamma_0)W(Z)\frac{\partial\rho(Y,Z;\gamma_0)}{\partial\gamma'}\right]. \tag{14}$$

*Furthermore,*

$$B = \operatorname*{plim}_{n\to\infty}\frac{1}{n}\sum_{i=1}^{n}\left[\frac{\partial\rho'(Y_i,Z_i;\hat\gamma_n)}{\partial\gamma}W(Z_i)\frac{\partial\rho(Y_i,Z_i;\hat\gamma_n)}{\partial\gamma'}\right] \tag{15}$$

*and*

$$4C = \operatorname*{plim}_{n\to\infty}\frac{1}{n}\frac{\partial Q_n(\hat\gamma_n)}{\partial\gamma}\frac{\partial Q_n(\hat\gamma_n)}{\partial\gamma'}, \tag{16}$$

*where*

$$\frac{\partial Q_n(\gamma)}{\partial\gamma} = 2\sum_{i=1}^{n}\frac{\partial\rho'(Y_i,Z_i;\gamma)}{\partial\gamma}W(Z_i)\rho(Y_i,Z_i;\gamma).$$

The MDE $\hat{\gamma}_n$ depends on the weight $W(Z)$. A natural question is how to choose $W(Z)$ to obtain the most efficient estimator. To answer this question, we first note that, since $\partial\rho'(Y,Z;\gamma_0)/\partial\gamma$ does not depend on $Y$, matrix $C$ in (14) can be written as

$$C = E\left[\frac{\partial\rho'(Y,Z;\gamma_0)}{\partial\gamma}W(Z)V(Z)W(Z)\frac{\partial\rho(Y,Z;\gamma_0)}{\partial\gamma'}\right],$$

where

$$V(Z) = E[\rho(Y,Z;\gamma_0)\rho'(Y,Z;\gamma_0)|Z]$$

and has elements

$$v_{11} = E[(Y - m_1(Z;\gamma_0))^2|Z],$$
$$v_{22} = E[(Y^2 - m_2(Z;\gamma_0))^2|Z]$$

and

$$v_{12} = E[(Y - m_1(Z;\gamma_0))(Y^2 - m_2(Z;\gamma_0))|Z].$$

Then, analogous to weighted (nonlinear) least squares estimation, we have

$$B^{-1}CB^{-1} \geq E\left[\frac{\partial\rho'(Y,Z;\gamma_0)}{\partial\gamma}V(Z)^{-1}\frac{\partial\rho(Y,Z;\gamma_0)}{\partial\gamma'}\right]^{-1} \tag{17}$$

(in the sense that the difference is nonnegative definite), and the lower bound is attained for $W(Z) = V(Z)^{-1}$ in $B$ and $C$. The matrix $V(Z)$ is invertible, if its determinant $v_{11}v_{22} - v_{12}^2 > 0$.

In general, $V(Z)$ is unknown, and it must be estimated before the MDE $\hat{\gamma}_n$ using $W(Z) = V(Z)^{-1}$ is computed. This can be done using the following two-stage procedure. First, minimize $Q_n(\gamma)$ using the identity matrix $W(Z) = I_2$ to obtain the first-stage estimator $\hat{\gamma}_n$. Second, estimate $V(Z)$ by

$$\hat{v}_{11} = \frac{1}{n}\sum_{i=1}^{n}(Y_i - m_1(Z_i;\hat{\gamma}_n))^2,$$

$$\hat{v}_{22} = \frac{1}{n}\sum_{i=1}^{n}(Y_i^2 - m_2(Z_i;\hat{\gamma}_n))^2$$

and

$$\hat{v}_{12} = \frac{1}{n}\sum_{i=1}^{n}(Y_i - m_1(Z_i;\hat{\gamma}_n))(Y_i^2 - m_2(Z_i;\hat{\gamma}_n)),$$

and then minimize $Q_n(\gamma)$ again with $W(Z) = \hat{V}(Z)^{-1}$ to obtain the two-stage estimator $\hat{\hat{\gamma}}_n$. Since the estimators $\hat{v}_{ij}$ are consistent for $v_{ij}$, the asymptotic covariance matrix of the two-stage estimator $\hat{\hat{\gamma}}_n$ is the same as the right-hand side of (17) and, therefore, $\hat{\hat{\gamma}}_n$ is asymptotically more efficient than the first-stage estimator $\hat{\gamma}_n$. More detailed discussions about the so-called feasible generalized least squares estimators can be found in, for example, Amemiya (1974) and Gallant [(1987), Chapter 5].

**4. Simulation-based estimator.** The MDE $\hat{\gamma}_n$ in the previous section is obtained by minimizing the objective function $Q_n(\gamma)$ in (13). The computation can be carried out using the usual numerical optimization procedures, if the explicit forms of $m_1(z;\gamma)$ and $m_2(z;\gamma)$ can be obtained. For some regression functions $g(x;\theta)$, however, explicit forms of the integrals in (11) and (12) may be difficult or impossible to derive. In this case, numerical integration techniques such as quadrature methods can be used. In practice, the numerical optimization of an objective function involving multiple integrals can be troublesome, especially when the dimension of the function is higher than three or four. To overcome this computational difficulty, in this section we consider a simulation-based approach for estimation in which the integrals are simulated by Monte Carlo methods such as importance sampling. This approach is similar to the method of simulated moments (MSM) of McFadden (1989) or Pakes and Pollard (1989).

The simulation-based estimator can be constructed in the following way. First, choose a known density function $\phi(t)$ and, for each $1 \le i \le n$, generate an i.i.d. random sample $\{t_{is}, s = 1, 2, \ldots, 2S\}$ from $\phi(t)$. Clearly, all samples $\{t_{is}, s = 1,2, \ldots, 2S, i = 1, 2, \ldots, n\}$ form a sequence of i.i.d. random variables. Then $m_1(z;\gamma)$ and $m_2(z;\gamma)$ can be approximated by the Monte Carlo simulators as

$$m_{1,S}(Z_i;\gamma) = \frac{1}{S}\sum_{s=1}^{S} \frac{g(Z_i + t_{is};\theta) f_\delta(t_{is};\psi)}{\phi(t_{is})},$$

$$m_{1,2S}(Z_i;\gamma) = \frac{1}{S}\sum_{s=S+1}^{2S} \frac{g(Z_i + t_{is};\theta) f_\delta(t_{is};\psi)}{\phi(t_{is})},$$

$$m_{2,S}(Z_i;\gamma) = \frac{1}{S}\sum_{s=1}^{S} \frac{g^2(Z_i + t_{is};\theta) f_\delta(t_{is};\psi)}{\phi(t_{is})} + \sigma_\varepsilon^2,$$

$$m_{2,2S}(Z_i;\gamma) = \frac{1}{S}\sum_{s=S+1}^{2S} \frac{g^2(Z_i + t_{is};\theta) f_\delta(t_{is};\psi)}{\phi(t_{is})} + \sigma_\varepsilon^2.$$

Therefore, a simulated version of the objective function $Q_n(\gamma)$ can be defined as

$$Q_{n,S}(\gamma) = \sum_{i=1}^{n} \rho_i^{(S)\prime}(\gamma) W(Z_i) \rho_i^{(2S)}(\gamma), \tag{18}$$

where

$$\rho_i^{(S)}(\gamma) = (Y_i - m_{1,S}(Z_i;\gamma), Y_i^2 - m_{2,S}(Z_i;\gamma))'$$

and

$$\rho_i^{(2S)}(\gamma) = (Y_i - m_{1,2S}(Z_i;\gamma), Y_i^2 - m_{2,2S}(Z_i;\gamma))'.$$

It is not difficult to see that $Q_{n,S}(\gamma)$ approximates $Q_n(\gamma)$ as $S \to \infty$, because by construction

$$E[m_{1,S}(Z_i;\gamma)|Z_i] = E[m_{1,2S}(Z_i;\gamma)|Z_i] = m_1(Z_i;\gamma)$$

and

$$E[m_{2,S}(Z_i;\gamma)|Z_i] = E[m_{2,2S}(Z_i;\gamma)|Z_i] = m_2(Z_i;\gamma).$$

In addition, $Q_{n,S}(\gamma)$ is an unbiased simulator for $Q_n(\gamma)$ in the sense that $EQ_{n,S}(\gamma) = EQ_n(\gamma)$, because, given $Y_i, Z_i$, $\rho_i^{(S)}(\gamma)$ and $\rho_i^{(2S)}(\gamma)$ are conditionally independent and hence

$$\begin{aligned} E[\rho_i^{(S)'}(\gamma)W(Z_i)\rho_i^{(2S)}(\gamma)] &= E[E(\rho_i^{(S)'}(\gamma)|Y_i,Z_i)W(Z_i)E(\rho_i^{(2S)}(\gamma)|Y_i,Z_i)] \\ &= E[\rho(Y_i,Z_i;\gamma)W(Z_i)\rho(Y_i,Z_i;\gamma)]. \end{aligned}$$

Finally, the simulation-based estimator (SE) for $\gamma$ is defined by

$$\hat{\gamma}_{n,S} = \underset{\gamma\in\Gamma}{\arg\min}\, Q_{n,S}(\gamma).$$

Note that, since $Q_{n,S}(\gamma)$ does not involve integrals any more, it is continuous in, and differentiable with respect to, $\gamma$, as long as functions $g(x;\theta)$ and $f_\delta(t;\psi)$ have these properties. In particular, the first derivative of $\rho_i^{(S)}(\gamma)$ becomes

$$\frac{\partial \rho_i^{(S)'}(\gamma)}{\partial\gamma} = -\left(\frac{\partial m_{1,S}(Z_i;\gamma)}{\partial\gamma}, \frac{\partial m_{2,S}(Z_i;\gamma)}{\partial\gamma}\right),$$

where $\partial m_{1,S}(Z_i;\gamma)/\partial\gamma$ is the column vector with elements

$$\begin{aligned} \frac{\partial m_{1,S}(Z_i;\gamma)}{\partial\theta} &= \frac{1}{S}\sum_{s=1}^{S}\frac{\partial g(Z_i+t_{is};\theta)}{\partial\theta}\frac{f_\delta(t_{is};\psi)}{\phi(t_{is})}, \\ \frac{\partial m_{1,S}(Z_i;\gamma)}{\partial\psi} &= \frac{1}{S}\sum_{s=1}^{S}\frac{g(Z_i+t_{is};\theta)}{\phi(t_{is})}\frac{\partial f_\delta(t_{is};\psi)}{\partial\psi}, \\ \frac{\partial m_{1,S}(Z_i;\gamma)}{\partial\sigma_\varepsilon^2} &= 0 \end{aligned}$$

and $\partial m_{2,S}(Z_i;\gamma)/\partial\gamma$ is the column vector with elements

$$\frac{\partial m_{2,S}(Z_i;\gamma)}{\partial\theta} = \frac{2}{S}\sum_{s=1}^{S}\frac{\partial g(Z_i+t_{is};\theta)}{\partial\theta}\frac{g(Z_i+t_{is};\theta)f_\delta(t_{is};\psi)}{\phi(t_{is})},$$

$$\frac{\partial m_{2,S}(Z_i;\gamma)}{\partial \psi} = \frac{1}{S}\sum_{s=1}^{S} \frac{g^2(Z_i+t_{is};\theta)}{\phi(t_{is})} \frac{\partial f_\delta(t_{is};\psi)}{\partial \psi},$$

$$\frac{\partial m_{2,S}(Z_i;\gamma)}{\partial \sigma_\varepsilon^2} = 1.$$

The derivatives $\partial m_{1,2S}(Z_i;\gamma)/\partial\gamma$ and $\partial m_{2,2S}(Z_i;\gamma)/\partial\gamma$ can be given similarly.

For the simulation-based estimator, we have the following results.

THEOREM 3. *Suppose the support of $\phi(t)$ covers the support of $f_\delta(t;\psi)$ for all $\psi \in \Psi$. Then the simulation estimator $\hat\gamma_{n,S}$ has the following properties:*

1. *Under Assumptions A1–A6, $\hat\gamma_{n,S} \xrightarrow{P} \gamma_0$ as $n \to \infty$.*
2. *Under Assumptions A1–A8, $\sqrt{n}(\hat\gamma_{n,S} - \gamma_0) \xrightarrow{L} N(0, B^{-1}C_S B^{-1})$, where*

$$(19)\qquad \begin{aligned} 2C_S &= E\left[\frac{\partial \rho_1^{(S)\prime}(\gamma_0)}{\partial\gamma} W(Z_1)\rho_1^{(2S)}(\gamma_0)\rho_1^{(2S)\prime}(\gamma_0)W(Z_1)\frac{\partial \rho_1^{(S)}(\gamma_0)}{\partial\gamma'}\right] \\ &\quad + E\left[\frac{\partial \rho_1^{(S)\prime}(\gamma_0)}{\partial\gamma} W(Z_1)\rho_1^{(2S)}(\gamma_0)\rho_1^{(S)\prime}(\gamma_0)W(Z_1)\frac{\partial \rho_1^{(2S)}(\gamma_0)}{\partial\gamma'}\right]. \end{aligned}$$

*Furthermore,*

$$(20)\qquad 4C_S = \operatorname*{plim}_{n\to\infty} \frac{1}{n} \frac{\partial Q_{n,S}(\hat\gamma_{n,S})}{\partial\gamma} \frac{\partial Q_{n,S}(\hat\gamma_{n,S})}{\partial\gamma'}.$$

In general, the simulation-based estimator $\hat\gamma_{n,S}$ is less efficient than the MDE $\hat\gamma_n$ of the previous section, due to the simulation approximation of $\rho_i(\gamma)$ by $\rho_i^{(S)}(\gamma)$ and $\rho_i^{(2S)}(\gamma)$. A natural question is how much efficiency is lost due to simulation. The following corollary provides an answer to this question.

COROLLARY 4. *Under the conditions of Theorem 3, it holds that*

$$(21)\qquad \begin{aligned} C_S &= C + \frac{1}{2S} E\left[\frac{\partial \rho_1' W(Z)(\rho_{11}-\rho_1)}{\partial\gamma} \frac{\partial (\rho_{11}-\rho_1)' W(Z)\rho_1}{\partial\gamma'}\right] \\ &\quad + \frac{1}{4S^2} E\left[\frac{\partial (\rho_{11}-\rho_1)' W(Z)(\rho_{12}-\rho_1)}{\partial\gamma} \frac{\partial (\rho_{12}-\rho_1)' W(Z)(\rho_{11}-\rho_1)}{\partial\gamma'}\right], \end{aligned}$$

*where $\rho_1 = \rho(Y_1, Z_1; \gamma_0)$ and*

$$\rho_{is} = \left(Y_i - \frac{g(Z_i+t_{is};\theta_0)f_\delta(t_{is};\psi_0)}{\phi(t_{is})}, Y_i^2 - \frac{g^2(Z_i+t_{is};\theta_0)f_\delta(t_{is};\psi_0)}{\phi(t_{is})} - \sigma_{\varepsilon 0}^2\right)'$$

*is the summand in $\rho_i^{(S)}(\gamma_0) = \sum_{s=1}^{S}\rho_{is}/S$.*

The above corollary shows that the efficiency loss caused by simulation has a magnitude of $O(1/S)$. Therefore, the larger the simulation size $S$, the smaller the efficiency loss. Furthermore, the efficiency loss reduces at rate $O(1/S)$. Asymptotically, the importance density $\phi(t)$ has no effect on the efficiency of the estimator, as long as it satisfies the condition of Theorem 3. In practice, however, the choice of $\phi(t)$ will affect the finite sample variances of the Monte Carlo estimators such as $m_{1,S}(Z_i;\gamma)$. Theoretically, the best choice of $\phi(t)$ is proportional to the absolute value of the integrand, which is $|g(z+t;\theta)f_\delta(t;\psi)|$ for $m_1(z;\gamma)$. Practically, however, a density close to being proportional to the integrand is a good choice. For more detailed discussion about importance sampling and variance reduction methods for numerical integration, see, for example, Evans and Swartz [(2000), Chapter 6].

In light of Corollary 4, the discussion in the previous section about the optimal choice of the weight $W(Z) = V(Z)^{-1}$ applies to the simulation-based estimator too, and will not be repeated here.

**5. Conclusion.** We have considered general nonlinear regression models with Berkson measurement errors in predictor variables. The measurement errors are assumed to have a general parametric distribution which is not necessarily normal, whereas the distribution of the random error in the regression equation is nonparametric. We have proposed a minimum distance estimator based on the first two conditional moments of the response variable given the observed predictor variables. We have shown that this estimator is consistent and asymptotically normally distributed under fairly general regularity conditions. To overcome the computational difficulty which may arise in the case where the objective function involves multiple integrals, a simulation-based estimator has been constructed. The consistency and asymptotic normality for this estimator have also been derived under regularity conditions similar to those for the minimum distance estimator. The results obtained generalize those of Wang (2003), which deals with the univariate model under normal distributions.

**6. Proofs.**

6.1. *Preliminary.* First, for ease of reading we restate some existing results which are used in the proofs. For this purpose, let $X = (X_1, X_2, \dots, X_n)$ be an i.i.d. random sample and let $\gamma$ be a vector of unknown parameters. Further, let $H(X_1, \gamma)$ and $S_n(X, \gamma)$ be measurable functions for any $\gamma \in \Gamma$, and be continuous in $\gamma \in \Gamma$ for almost all possible values of $X$. In addition, the parameter space $\Gamma \subset \mathbb{R}^d$ is compact. Using this notation, Theorems 4.2.1, 4.1.1 and 4.1.5 of Amemiya (1985) can be stated as follows.

LEMMA 5. *Suppose* $E\sup_{\gamma\in\Gamma}|H(X_1,\gamma)| < \infty$. *Then* $\frac{1}{n}\sum_{i=1}^n H(X_i,\gamma)$ *converges in probability to* $EH(X_1,\gamma)$ *uniformly in* $\gamma \in \Gamma$.

LEMMA 6. *Suppose, as* $n \to \infty$, $S_n(X,\gamma)$ *converges in probability to a nonstochastic function* $S(\gamma)$ *uniformly in* $\gamma \in \Gamma$, *and* $S(\gamma)$ *attains a unique minimum at* $\gamma_0 \in \Gamma$. *Then the estimator* $\hat{\gamma}_n = \arg\min_{\gamma\in\Gamma} S_n(X,\gamma)$ *converges in probability to* $\gamma_0$.

LEMMA 7. *Suppose, as* $n \to \infty$, $S_n(X,\gamma)$ *converges in probability to a nonstochastic function* $S(\gamma)$ *uniformly in* $\gamma$ *in an open neighborhood of* $\gamma_0$, *and* $S(\gamma)$ *is continuous at* $\gamma_0$. *Then* $\operatorname{plim}_{n\to\infty} \hat{\gamma}_n = \gamma_0$ *implies* $\operatorname{plim}_{n\to\infty} S_n(X,\hat{\gamma}_n) = S(\gamma_0)$.

To simplify the notation in the proofs, we will denote $\rho(Y_i, Z_i; \gamma)$ as $\rho_i(\gamma)$, and $W(Z_i)$ as $W_i$, as far as these cause no confusion. For any matrix $A$, its Euclidean norm is denoted as $\|A\| = \sqrt{\operatorname{trace}(A'A)}$, and $\operatorname{vec}A$ denotes the column vector consisting of the columns of $A$. Further, $\otimes$ denotes the Kronecker product operator.

PROOF OF THEOREM 1. We show that Assumptions A1–A6 are sufficient for all conditions of Lemma 6. First, by Hölder's inequality and Assumption A5 we have

$$E \int \sup_{\Theta\times\Psi} |g(Z+t;\theta)|^j f_\delta(t;\psi)\,dt < \infty \tag{22}$$

for $j = 1, 2, 3$. It follows from Assumptions A1, A2 and the dominated convergence theorem that $m_1(z;\gamma)$, $m_2(z;\gamma)$ and therefore $Q_n(\gamma)$ are continuous in $\gamma \in \Gamma$. Let

$$Q(\gamma) = E\rho_1'(\gamma) W(Z_1)\rho_1(\gamma).$$

Again by Hölder's inequality, (22) and Assumption A3 we have

$$\begin{aligned}
&E\|W_1\| \sup_{\Gamma} [Y_1 - m_1(Z_1;\gamma)]^2 \\
&\quad \le 2E\|W_1\|Y_1^2 + 2E\|W_1\| \sup_{\Gamma} m_1^2(Z_1;\gamma) \\
&\quad \le 2E\|W_1\|Y_1^2 + 2E\|W_1\| \int \sup_{\Theta\times\Psi} g^2(Z_1+t;\theta) f_\delta(t;\psi)\,dt \\
&\quad < \infty
\end{aligned}$$

and

$$\begin{aligned}
&E\|W_1\| \sup_{\Gamma} [Y_1^2 - m_2(Z_1;\gamma)]^2 \\
&\quad \le 3E\|W_1\|Y_1^4 + 3E\|W_1\| \int \sup_{\Theta\times\Psi} g^4(Z_1+t;\theta) f_\delta(t;\psi)\,dt \\
&\qquad + 3E\|W_1\| \sup_{\Sigma} \sigma_\varepsilon^4 < \infty,
\end{aligned}$$

which imply

$$E\sup_{\Gamma}\rho_1'(\gamma)W_1\rho_1(\gamma)\le E\|W_1\|\sup_{\Gamma}\|\rho_1(\gamma)\|^2<\infty.$$

It follows from Lemma 5 that $\frac{1}{n}Q_n(\gamma)$ converges in probability to $Q(\gamma)$ uniformly in $\gamma\in\Gamma$. Further, since

$$\begin{aligned}E[\rho_1'(\gamma_0)W_1(\rho_1(\gamma)-\rho_1(\gamma_0))]&=E[E(\rho_1'(\gamma_0)|Z_1)W_1(\rho_1(\gamma)-\rho_1(\gamma_0))]\\&=0,\end{aligned}$$

we have

$$Q(\gamma)=Q(\gamma_0)+E[(\rho_1(\gamma)-\rho_1(\gamma_0))'W_1(\rho_1(\gamma)-\rho_1(\gamma_0))].$$

It follows that $Q(\gamma)\ge Q(\gamma_0)$ and, by Assumption A6, equality holds if and only if $\gamma=\gamma_0$. Thus all conditions of Lemma 6 hold and, therefore, $\hat\gamma_n\xrightarrow{P}\gamma_0$ follows. □

Proof of Theorem 2. By Assumption A7 the first derivative $\partial Q_n(\gamma)/\partial\gamma$ exists and has a first-order Taylor expansion in a neighborhood $\Gamma_0\subset\Gamma$ of $\gamma_0$. Since $\partial Q_n(\hat\gamma_n)/\partial\gamma=0$ and $\hat\gamma_n\xrightarrow{P}\gamma_0$, for sufficiently large $n$ we have

$$0=\frac{\partial Q_n(\gamma_0)}{\partial\gamma}+\frac{\partial^2Q_n(\tilde\gamma_n)}{\partial\gamma\,\partial\gamma'}(\hat\gamma_n-\gamma_0),\tag{23}$$

where $\|\tilde\gamma_n-\gamma_0\|\le\|\hat\gamma_n-\gamma_0\|$. The first derivative of $Q_n(\gamma)$ in (23) is given by

$$\frac{\partial Q_n(\gamma)}{\partial\gamma}=2\sum_{i=1}^n\frac{\partial\rho_i'(\gamma)}{\partial\gamma}W_i\rho_i(\gamma),$$

where

$$\frac{\partial\rho_i'(\gamma)}{\partial\gamma}=-\left(\frac{\partial m_1(Z_i;\gamma)}{\partial\gamma},\frac{\partial m_2(Z_i;\gamma)}{\partial\gamma}\right)$$

and the first derivatives of $m_1(Z_i;\gamma)$ and $m_2(Z_i;\gamma)$ with respect to $\gamma$ are given in Assumption A8. Therefore, by the central limit theorem we have

$$\frac{1}{\sqrt n}\frac{\partial Q_n(\gamma_0)}{\partial\gamma}\xrightarrow{L}N(0,4C),\tag{24}$$

where

$$C=E\left[\frac{\partial\rho_i'(\gamma_0)}{\partial\gamma}W_i\rho_i(\gamma_0)\rho_i'(\gamma_0)W_i\frac{\partial\rho_i(\gamma_0)}{\partial\gamma'}\right],$$

as is given in (14). The second derivative of $Q_n(\gamma)$ in (23) is given by

$$\frac{\partial^2Q_n(\gamma)}{\partial\gamma\,\partial\gamma'}=2\sum_{i=1}^n\left[\frac{\partial\rho_i'(\gamma)}{\partial\gamma}W_i\frac{\partial\rho_i(\gamma)}{\partial\gamma'}+(\rho_i'(\gamma)W_i\otimes I_{p+q+1})\frac{\partial\operatorname{vec}(\partial\rho_i'(\gamma)/\partial\gamma)}{\partial\gamma'}\right],$$

where

$$\frac{\partial\,\mathrm{vec}(\partial\rho_i'(\gamma)/\partial\gamma)}{\partial\gamma'} = -\left(\frac{\partial^2 m_1(z;\gamma)}{\partial\gamma\,\partial\gamma'}, \frac{\partial^2 m_2(z;\gamma)}{\partial\gamma\,\partial\gamma'}\right)'.$$

Again, by Assumption A7, the nonzero elements in $\partial^2 m_1(z;\gamma)/\partial\gamma\,\partial\gamma'$ are

$$\begin{aligned}
\frac{\partial^2 m_1(z;\gamma)}{\partial\theta\,\partial\theta'} &= \int \frac{\partial^2 g(z+t;\theta)}{\partial\theta\,\partial\theta'} f_\delta(t;\psi)\,dt,\\
\frac{\partial^2 m_1(z;\gamma)}{\partial\theta\,\partial\psi'} &= \int \frac{\partial g(z+t;\theta)}{\partial\theta}\frac{\partial f_\delta(t;\psi)}{\partial\psi'}\,dt,\\
\frac{\partial^2 m_1(z;\gamma)}{\partial\psi\,\partial\psi'} &= \int g(z+t;\theta)\frac{\partial^2 f_\delta(t;\psi)}{\partial\psi\,\partial\psi'}\,dt,
\end{aligned}$$

and the nonzero elements in $\partial^2 m_2(z;\gamma)/\partial\gamma\,\partial\gamma'$ are

$$\begin{aligned}
\frac{\partial^2 m_2(z;\gamma)}{\partial\theta\,\partial\theta'} &= 2\int \frac{\partial^2 g(z+t;\theta)}{\partial\theta\,\partial\theta'} g(z+t;\theta) f_\delta(t;\psi)\,dt\\
&\quad + 2\int \frac{\partial g(z+t;\theta)}{\partial\theta}\frac{\partial g(z+t;\theta)}{\partial\theta'} f_\delta(t;\psi)\,dt,\\
\frac{\partial^2 m_2(z;\gamma)}{\partial\theta\,\partial\psi'} &= 2\int g(z+t;\theta)\frac{\partial g(z+t;\theta)}{\partial\theta}\frac{\partial f_\delta(t;\psi)}{\partial\psi'}\,dt,\\
\frac{\partial^2 m_2(z;\gamma)}{\partial\psi\,\partial\psi'} &= \int g^2(z+t;\theta)\frac{\partial^2 f_\delta(t;\psi)}{\partial\psi\,\partial\psi'}\,dt.
\end{aligned}$$

Analogously to the proof of Theorem 1, we can verify by Assumption A7 and Lemma 5 that $(1/n)\,\partial^2 Q_n(\gamma)/\partial\gamma\,\partial\gamma'$ converges in probability to $\partial^2 Q(\gamma)/\partial\gamma\,\partial\gamma'$ uniformly in $\gamma\in\Gamma_0$. Therefore, by Lemma 7 we have

$$\begin{aligned}
&\frac{1}{n}\frac{\partial^2 Q_n(\tilde\gamma_n)}{\partial\gamma\,\partial\gamma'}\\
&\xrightarrow{P} 2E\left[\frac{\partial\rho_1'(\gamma_0)}{\partial\gamma}W_1\frac{\partial\rho_1(\gamma_0)}{\partial\gamma'} + (\rho_1'(\gamma_0)W_1\otimes I_{p+q+1})\frac{\partial\,\mathrm{vec}(\partial\rho_1'(\gamma_0)/\partial\gamma)}{\partial\gamma'}\right]\\
&= 2B,
\end{aligned}\tag{25}$$

where the second equality holds, because

$$\begin{aligned}
&E\left[(\rho_1'(\gamma_0)W_1\otimes I_{p+q+1})\frac{\partial\,\mathrm{vec}(\partial\rho_1'(\gamma_0)/\partial\gamma)}{\partial\gamma'}\right]\\
&\quad = E\left[(E(\rho_1'(\gamma_0)|Z_1)W_1\otimes I_{p+q+1})\frac{\partial\,\mathrm{vec}(\partial\rho_1'(\gamma_0)/\partial\gamma)}{\partial\gamma'}\right]\\
&\quad = 0.
\end{aligned}$$

It follows then from (23)–(25), Assumption A8 and the Slutsky theorem that $\sqrt{n}(\hat{\gamma}_n - \gamma_0) \xrightarrow{L} N(0, B^{-1}CB^{-1})$. Finally, (15) and (16) can be similarly verified by Lemma 7. □

Proof of Theorem 3. The proof for part 1 of Theorem 3 is analogous to that for Theorem 1. First, Assumptions A1 and A2 imply that $Q_{n,S}(\gamma)$ is continuous in $\gamma \in \Gamma$. Then, by Lemma 5 we have, as $n \to \infty$, uniformly in $\gamma \in \Gamma$ that

$$\begin{aligned}\frac{1}{n} Q_{n,S}(\gamma) \xrightarrow{P} E[\rho_1^{(S)'}(\gamma) W(Z_1) \rho_1^{(2S)}(\gamma)]& \\ &= E[E(\rho_1^{(S)'}(\gamma)|Y_1, Z_1) W(Z_1) E(\rho_1^{(2S)}(\gamma)|Y_1, Z_1)] \\ &= E[\rho_1'(\gamma) W(Z_1) \rho_1(\gamma)] \\ &= Q(\gamma).\end{aligned}$$

Finally, $\hat{\gamma}_{n,S} \xrightarrow{P} \gamma_0$ follows from Assumption A6 and Lemma 6.

The proof of part 2 of Theorem 3 is analogous to that of Theorem 2. First, by Assumption A7 we have the first-order Taylor expansion of $\partial Q_{n,S}(\gamma)/\partial\gamma$ in a neighborhood $\Gamma_0 \subset \Gamma$ of $\gamma_0$,

$$0 = \frac{\partial Q_{n,S}(\gamma_0)}{\partial \gamma} + \frac{\partial^2 Q_{n,S}(\tilde{\gamma}_{n,S})}{\partial \gamma \, \partial \gamma'} (\hat{\gamma}_{n,S} - \gamma_0), \tag{26}$$

where $\|\tilde{\gamma}_{n,S} - \gamma_0\| \le \|\hat{\gamma}_{n,S} - \gamma_0\|$ and the first derivative of $Q_{n,S}(\gamma)$ is given by

$$\frac{\partial Q_{n,S}(\gamma)}{\partial \gamma} = \sum_{i=1}^{n} \left[ \frac{\partial \rho_i^{(S)'}(\gamma)}{\partial \gamma} W_i \rho_i^{(2S)}(\gamma) + \frac{\partial \rho_i^{(2S)'}(\gamma)}{\partial \gamma} W_i \rho_i^{(S)}(\gamma) \right].$$

Since $\rho_i^{(S)}(\gamma)$ has the same distribution as $\rho_i^{(2S)}(\gamma)$, all terms in the above summation are i.i.d. and have the common covariance matrix

$$\begin{aligned}&E\left[ \frac{\partial \rho_i^{(S)'}(\gamma_0)}{\partial \gamma} W_i \rho_i^{(2S)}(\gamma_0) \rho_i^{(2S)'}(\gamma_0) W_i \frac{\partial \rho_i^{(S)}(\gamma_0)}{\partial \gamma'} \right] \\ &\quad + E\left[ \frac{\partial \rho_i^{(S)'}(\gamma_0)}{\partial \gamma} W_i \rho_i^{(2S)}(\gamma_0) \rho_i^{(S)'}(\gamma_0) W_i \frac{\partial \rho_i^{(2S)}(\gamma_0)}{\partial \gamma'} \right] \\ &\quad + E\left[ \frac{\partial \rho_i^{(2S)'}(\gamma_0)}{\partial \gamma} W_i \rho_i^{(S)}(\gamma_0) \rho_i^{(2S)'}(\gamma_0) W_i \frac{\partial \rho_i^{(S)}(\gamma_0)}{\partial \gamma'} \right] \\ &\quad + E\left[ \frac{\partial \rho_i^{(2S)'}(\gamma_0)}{\partial \gamma} W_i \rho_i^{(S)}(\gamma_0) \rho_i^{(S)'}(\gamma_0) W_i \frac{\partial \rho_i^{(2S)}(\gamma_0)}{\partial \gamma'} \right] \\ &\quad = 4 C_S.\end{aligned}$$

It follows by the central limit theorem that, as $n \to \infty$,

$$\frac{1}{\sqrt{n}} \frac{\partial Q_{n,S}(\gamma_0)}{\partial \gamma} \xrightarrow{L} N(0, 4C_S). \tag{27}$$

Now, the second derivative in (26) is given by

$$\begin{aligned}
\frac{\partial^2 Q_{n,S}(\gamma)}{\partial \gamma\, \partial \gamma'} = \sum_{i=1}^{n} \Big[ & \frac{\partial \rho_i^{(S)'}(\gamma)}{\partial \gamma} W_i \frac{\partial \rho_i^{(2S)}(\gamma)}{\partial \gamma'} \\
& + (\rho_i^{(2S)'}(\gamma) W_i \otimes I_{p+q+1}) \frac{\partial \operatorname{vec}(\partial \rho_i^{(S)'}(\gamma)/\partial \gamma)}{\partial \gamma'} \Big] \\
+ \sum_{i=1}^{n} \Big[ & \frac{\partial \rho_i^{(2S)'}(\gamma)}{\partial \gamma} W_i \frac{\partial \rho_i^{(S)}(\gamma)}{\partial \gamma'} \\
& + (\rho_i^{(S)'}(\gamma) W_i \otimes I_{p+q+1}) \frac{\partial \operatorname{vec}(\partial \rho_i^{(2S)'}(\gamma)/\partial \gamma)}{\partial \gamma'} \Big],
\end{aligned}$$

where

$$\frac{\partial \operatorname{vec}(\partial \rho_i^{(S)'}(\gamma)/\partial \gamma)}{\partial \gamma'} = -\left( \frac{\partial^2 m_{1,S}(z;\gamma)}{\partial \gamma\, \partial \gamma'}, \frac{\partial^2 m_{2,S}(z;\gamma)}{\partial \gamma\, \partial \gamma'} \right)'$$

and the nonzero elements in $\partial^2 m_{1,S}(Z_i;\gamma)/\partial \gamma\, \partial \gamma'$ are

$$\frac{\partial^2 m_{1,S}(Z_i;\gamma)}{\partial \theta\, \partial \theta'} = \frac{1}{S} \sum_{s=1}^{S} \frac{\partial^2 g(Z_i + t_{is};\theta)}{\partial \theta\, \partial \theta'} \frac{f_\delta(t_{is};\psi)}{\phi(t_{is})},$$

$$\frac{\partial^2 m_{1,S}(Z_i;\gamma)}{\partial \theta\, \partial \psi'} = \frac{1}{S} \sum_{s=1}^{S} \frac{\partial g(Z_i + t_{is};\theta)}{\partial \theta} \frac{\partial f_\delta(t_{is};\psi)}{\partial \psi'} \frac{1}{\phi(t_{is})},$$

$$\frac{\partial^2 m_{1,S}(Z_i;\gamma)}{\partial \psi\, \partial \psi'} = \frac{1}{S} \sum_{s=1}^{S} \frac{g(Z_i + t_{is};\theta)}{\phi(t_{is})} \frac{\partial^2 f_\delta(t_{is};\psi)}{\partial \psi\, \partial \psi'},$$

and the nonzero elements in $\partial^2 m_{2,S}(Z_i;\gamma)/\partial \gamma \partial \gamma'$ are

$$\begin{aligned}
\frac{\partial^2 m_{2,S}(Z_i;\gamma)}{\partial \theta\, \partial \theta'} = \frac{2}{S} \sum_{s=1}^{S} \Big[ & \frac{\partial^2 g(Z_i + t_{is};\theta)}{\partial \theta\, \partial \theta'} \frac{g(Z_i + t_{is};\theta) f_\delta(t_{is};\psi)}{\phi(t_{is})} \\
& + \frac{\partial g(Z_i + t_{is};\theta)}{\partial \theta} \frac{\partial g(Z_i + t_{is};\theta)}{\partial \theta'} \frac{f_\delta(t_{is};\psi)}{\phi(t_{is})} \Big],
\end{aligned}$$

$$\frac{\partial^2 m_{2,S}(Z_i;\gamma)}{\partial \theta\, \partial \psi'} = \frac{2}{S} \sum_{s=1}^{S} \frac{g(Z_i + t_{is};\theta)}{\phi(t_{is})} \frac{\partial g(Z_i + t_{is};\theta)}{\partial \theta} \frac{\partial f_\delta(t_{is};\psi)}{\partial \psi'},$$

$$\frac{\partial^2 m_{2,S}(Z_i;\gamma)}{\partial \psi\, \partial \psi'} = \frac{1}{S} \sum_{s=1}^{S} \frac{g^2(Z_i + t_{is};\theta)}{\phi(t_{is})} \frac{\partial^2 f_\delta(t_{is};\psi)}{\partial \psi\, \partial \psi'}.$$

Again, by Assumption A7 and Lemma 7, uniformly in $\gamma \in \Gamma$,

$$\begin{aligned}\frac{\partial^2 Q_{n,S}(\tilde\gamma_n)}{\partial\gamma\,\partial\gamma'} \xrightarrow{P} E\Big[&\frac{\partial\rho_1^{(S)'}(\gamma_0)}{\partial\gamma}W_1\frac{\partial\rho_1^{(2S)}(\gamma_0)}{\partial\gamma'} \\ &+ (\rho_1^{(2S)'}(\gamma_0)W_1\otimes I_{p+q+1})\frac{\partial\,\mathrm{vec}(\partial\rho_1^{(S)'}(\gamma_0)/\partial\gamma)}{\partial\gamma'}\Big] \\ + E\Big[&\frac{\partial\rho_1^{(2S)'}(\gamma_0)}{\partial\gamma}W_1\frac{\partial\rho_1^{(S)}(\gamma_0)}{\partial\gamma'} \\ &+ (\rho_1^{(S)'}(\gamma_0)W_1\otimes I_{p+q+1})\frac{\partial\,\mathrm{vec}(\partial\rho_1^{(2S)'}(\gamma_0)/\partial\gamma)}{\partial\gamma'}\Big] \\ = E\Big[&\frac{\partial\rho_1^{(S)'}(\gamma_0)}{\partial\gamma}W_1\frac{\partial\rho_1^{(2S)}(\gamma_0)}{\partial\gamma'} + \frac{\partial\rho_1^{(2S)'}(\gamma_0)}{\partial\gamma}W_1\frac{\partial\rho_1^{(S)}(\gamma_0)}{\partial\gamma'}\Big] \\ = 2E\Big[&\frac{\partial\rho_1^{(S)'}(\gamma_0)}{\partial\gamma}W_1\frac{\partial\rho_1^{(2S)}(\gamma_0)}{\partial\gamma'}\Big] \\ = 2B&,\end{aligned} \tag{28}$$

where the first equality follows from

$$\begin{aligned}E\Big[&(\rho_1^{(2S)'}(\gamma_0)W_1\otimes I_{p+q+1})\frac{\partial\,\mathrm{vec}(\partial\rho_1^{(S)'}(\gamma_0)/\partial\gamma)}{\partial\gamma'}\Big] \\ &= E\Big[E(\rho_1^{(2S)'}(\gamma)|Z_i)W_1\otimes I_{p+q+1}\frac{\partial\,\mathrm{vec}(\partial\rho_1^{(S)'}(\gamma)/\partial\gamma)}{\partial\gamma'}\Big] \\ &= 0,\end{aligned}$$

and the last equality holds because

$$\begin{aligned}E\Big[\frac{\partial\rho_1^{(S)'}(\gamma_0)}{\partial\gamma}W_1\frac{\partial\rho_1^{(2S)}(\gamma_0)}{\partial\gamma'}\Big] &= E\Big[E\Big(\frac{\partial\rho_1^{(S)'}(\gamma_0)}{\partial\gamma}\Big|Z_1\Big)W_1E\Big(\frac{\partial\rho_1^{(2S)}(\gamma_0)}{\partial\gamma'}\Big|Z_1\Big)\Big] \\ &= E\Big[\frac{\partial\rho_1'(\gamma_0)}{\partial\gamma}W_1\frac{\partial\rho_1(\gamma_0)}{\partial\gamma'}\Big].\end{aligned}$$

By (26)–(28) and the Slutsky theorem, we have $\sqrt{n}(\hat\gamma_{n,S}-\gamma_0) \xrightarrow{L} N(0, B^{-1}C_SB^{-1})$. Finally, (20) can be similarly shown by Lemma 7. $\square$

PROOF OF COROLLARY 4. To simplify notation, in the following we denote $\rho_i = \rho_i(\gamma_0)$ and, correspondingly, $\rho_i^{(S)} = \rho_i^{(S)}(\gamma_0)$. Then the common term of $\partial Q_{n,S}(\gamma_0)/\partial\gamma$ in (26) can be written as

$$T = \frac{\partial\rho_i^{(S)'}W_i\rho_i^{(2S)}}{\partial\gamma}$$

$$= T_1 + T_2 + T_3,$$

where

$$T_1 = \frac{\partial \rho_i' W_i \rho_i}{\partial \gamma},$$

$$T_2 = \frac{\partial \rho_i' W_i(\rho_i^{(S)} - \rho_i)}{\partial \gamma} + \frac{\partial \rho_i' W_i(\rho_i^{(2S)} - \rho_i)}{\partial \gamma}$$

and

$$T_3 = \frac{\partial (\rho_i^{(S)} - \rho_i)' W_i(\rho_i^{(2S)} - \rho_i)}{\partial \gamma}.$$

Since $\rho_i^{(S)}$ and $\rho_i^{(2S)}$ are conditionally independent given $Y_i$ and $Z_i$, $T_1$, $T_2$ and $T_3$ are mutually uncorrelated and hence

$$E(TT') = E(T_1T_1') + E(T_2T_2') + E(T_3T_3'), \tag{29}$$

where $E(TT') = 4C_S$ and $E(T_1T_1') = 4C$. Furthermore, since $\rho_i^{(S)}$ and $\rho_i^{(2S)}$ have the same distribution,

$$E(T_2T_2') = 2E\left[\frac{\partial \rho_i' W_i(\rho_i^{(S)} - \rho_i)}{\partial \gamma} \frac{\partial (\rho_i^{(S)} - \rho_i)' W_i \rho_i}{\partial \gamma'}\right].$$

Now write $\rho_i^{(S)}(\gamma_0) = \frac{1}{S}\sum_{s=1}^S \rho_{is}$, where

$$\rho_{is} = \left(Y_i - \frac{g(Z_i + t_{is}; \theta_0) f_\delta(t_{is}; \psi_0)}{\phi(t_{is})}, Y_i^2 - \frac{g^2(Z_i + t_{is}; \theta_0) f_\delta(t_{is}; \psi_0)}{\phi(t_{is})} - \sigma_{\varepsilon 0}^2\right)'.$$

Then, since $\rho_{is}, s = 1, 2, \ldots, S$, are independent given $Y_i, Z_i$, we have

$$\begin{aligned} E(T_2T_2') &= \frac{2}{S^2} E\left[\sum_{s=1}^S \frac{\partial \rho_i' W_i(\rho_{is} - \rho_i)}{\partial \gamma} \sum_{s=1}^S \frac{\partial (\rho_{is} - \rho_i)' W_i \rho_i}{\partial \gamma'}\right] \\ &= \frac{2}{S} E\left[\frac{\partial \rho_1' W_1(\rho_{11} - \rho_1)}{\partial \gamma} \frac{\partial (\rho_{11} - \rho_1)' W_1 \rho_1}{\partial \gamma'}\right]. \end{aligned} \tag{30}$$

In the same way, we can show that

$$\begin{aligned} &E(T_3T_3') \\ &\quad = \frac{1}{S^2} E\left[\frac{\partial (\rho_{11} - \rho_1)' W(Z)(\rho_{12} - \rho_1)}{\partial \gamma} \frac{\partial (\rho_{12} - \rho_1)' W(Z)(\rho_{11} - \rho_1)}{\partial \gamma'}\right]. \end{aligned} \tag{31}$$

The corollary follows from (29)–(31). □

**Acknowledgments.** The author is grateful to the Editor, John I. Marden, an Associate Editor and a referee for their constructive comments and suggestions which led to the improvement of an earlier version of this paper.

Department of Statistics
University of Manitoba
Winnipeg, Manitoba
Canada R3T 2N2
e-mail: liqun_wang@umanitoba.ca